\documentclass[a4paper,12pt,leqno]{article}
\usepackage{amsmath,amsthm,amsfonts,amssymb}
\usepackage[titletoc,toc,title]{appendix}
\numberwithin{equation}{section}

\newcommand{\calH}{\mathcal{H}}

\newcommand{\Comp}{\mathbb{C}}

\DeclareMathOperator{\diag}{diag}
\DeclareMathOperator{\End}{End}
\renewcommand{\epsilon}{\varepsilon}
\DeclareMathOperator{\Hom}{Hom}
\DeclareMathOperator{\Id}{Id}
\renewcommand{\Im}{\mathrm{Im}}
\newcommand{\Integer}{\mathbb{Z}}

\renewcommand{\Re}{\mathrm{Re}}
\newcommand{\Real}{\mathbb{R}}
\newcommand{\tensor}{\otimes}
\DeclareMathOperator{\tr}{tr}

\newcommand\lemref[1]{Lemma~\ref{#1}}

\newcommand\remref[1]{Remark~\ref{#1}}
\newcommand\secref[1]{\S\ref{#1}}
\newcommand\appref[1]{Appendix~\ref{#1}}
\newcommand\tableref[1]{Table~\ref{#1}}

\newtheorem{theorem}{Theorem}[section]

\newtheorem{lemma}[theorem]{Lemma}

\theoremstyle{definition}

\theoremstyle{remark} \newtheorem{remark}[theorem]{Remark}
\newcommand\comment[1]{}
\newcommand\commentout[1]{}
\begin{document}
\title%
{$Q$-operators for higher spin \\eight vertex models\\
with a rational anisotropy parameter}
\author{Takashi TAKEBE\thanks{e-mail: \tt ttakebe@hse.ru}\\
Faculty of Mathematics,\\
International Laboratory of Representation Theory\\ 
and Mathematical Physics\\
National Research University Higher School of Economics\\
Usacheva Street 6, Moscow, 119048, Russia.}
\date{}



\maketitle
\begin{abstract}
 $Q$-operators for generalised eight vertex models associated to higher
 spin representations of the Sklyanin algebra are constructed by
 Baxter's first method and Fabricius's method, when the anisotropy
 parameter is rational.
\end{abstract}

\section{Introduction}
\label{sec:intro}

This paper is a continuation of \cite{tak:16}, in which we constructed a
$Q$-operator of the highger spin generalisation of the eight vertex
model. In the present paper we construct the $Q$-operator of the same
model with a rational anisotropy parameter $\eta$, following Baxter's
1972 paper \cite{bax:72} and Fabricius's paper \cite{fab:07}.

Baxter \cite{bax:72} introduced the $Q$-operator as an auxiliary
tool to find eigenvalues of the transfer matrix of the eight vertex
model. Its fundamental property is the {\em $TQ$-relation},
\begin{equation*}
    T(u) Q(u) = Q(u) T(u) = h_-(u) Q(u-2\eta) + h_+(u) Q(u+2\eta),
\end{equation*}
with the transfer matrix $T(u)$, the commutativity with itself,
$[Q(u),Q(u')]=0$, and the quasi-periodicity with respect to the shift of
the spectral parameter $u\mapsto u+1$, $u\mapsto u+\tau$.

We constructed such an operator for the higher spin eight vertex models
in \cite{tak:16}, following Baxter's 1973 paper \cite{bax:73} and his
book \cite{bax:82}, when the number of sites $N$ was even. This model is
defined by spin $l$ representations of the Sklyanin algebra
\cite{skl:82}, \cite{skl:83} and has been studied by the author in
\cite{tak:92}, \cite{tak:95}, \cite{tak:96} by means of the
generalised algebraic Bethe Ansatz of Takhtajan and Faddeev
\cite{takh-fad:79}. The results of \cite{tak:16} not only confirmed the
results obtained by the algebraic Bethe Ansatz in \cite{tak:92}, but
also proved the sum rule of Behte roots conjectured in \cite{tak:95}.

In the present paper we show that Baxter's very first construction of
the $Q$-operator in \cite{bax:72} and Fabricius's new construction in
\cite{fab:07} (see also \cite{roa:07-2}) can be applied to the higher
spin model with a rational anisotropy parameter.  Recall that in
\cite{tak:92} and \cite{tak:95} the generalised Bethe Ansatz was shown
to be applicable to the higher spin model under the assumption,
$Nl\in\Integer$, but our construction of the $Q$-operator in
\cite{tak:16} was valid only for even $N$. In particular, the case when
the number of sites $N$ is odd and the spin $l$ is an integer was
excluded in \cite{tak:16}, while the generalised Bethe Ansatz was
applied in \cite{tak:92} and \cite{tak:95} to such a case as well. We
show that Baxter's construction in \cite{bax:72} gives the $Q$-operator
also in this case. Thus we fill the gap between the generalised Bethe
Ansatz and the $Q$-operator approach for the higher spin model.

As was pointed out by Fabricius and McCoy in \cite{fab-mcc:03-05},
Baxter's 1972 construction has a weak point. In \cite{bax:72} the
$Q$-operator $Q(u)$ was defined by $Q(u)=Q_R(u)\, Q_R(u_0)^{-1}$, where
$Q_R(u)$ is an intermediate object in the construction and $u_0$ is a
point, at which $Q_R(u)$ is invertible. Fabricius and McCoy found
numerically that $Q_R(u)$ might be degenerate for certain $\eta$, which
means that this construction cannot be applied.

We do not touch upon this problem and assume that $Q_R(u)$ is
non-degenerate at generic $u$. The reason is not because we think this
problem negligible, but, on the contrary, because it requires serious
study, to which another article should be devoted. Anyway, according to
\cite{fab-mcc:07}, either Baxter's or Fabricius's construction works for
any rational anisotropy parameter\footnote{The normalisation of theta
functions by Fabricius and McCoy in \cite{fab-mcc:03-05},
\cite{fab-mcc:07}, \cite{fab-mcc:09} and that by Sklyanin in
\cite{skl:83} differs by a modular transformation. (See the explicit
form of the $R$-matrix, \eqref{R} in this paper.) Since we use
Sklyanin's normalisation, our cases correspond to $m_1=0$ in
\cite{fab-mcc:03-05}.}.

\medskip
What is remarkable in our construction is that Baxter's and Fabricius's
methods work almost just as they are, although they seem to make use of
special form of matrix elements.

\medskip
This paper is organised as follows: in \secref{sec:model} we recall the
model introduced in \cite{tak:92} by defining its transfer matrix and
general strategy of the construction of the $Q$-operator, which is
essentially due to Baxter \cite{bax:72}: First we construct an operator
$Q_R(u)$ which satisfies the $TQ$-relation \eqref{TQ} as the trace of a
product of auxiliary operators $S(u)$. Then by transposing $Q_R(u)$ (in
our case, by taking the Hermitian adjoint with respect to the Sklyanin
form) we obtain an operator $Q_L(u)$ satisfying the $QT$-relation
\eqref{QT}. Operators $Q_R(u)$ and $Q_L(u)$ satisfy the commutation
relation \eqref{QLQR=QLQR}, from which follows the definition of the
$Q$-operator, \eqref{Q=QRQR-1=QL-1QL}.

In \secref{sec:baxter72} we construct a $Q$-operator, following Baxter's
1972 paper \cite{bax:72}, and in \secref{sec:fabricius} we follow
essentially Fabricius's paper \cite{fab:07}. The auxiliary
operator-valued matrix $S(u)$ in \secref{sec:baxter72} is tridiagonal in
ordinary sense \eqref{S:matrix:bax} and $S(u)$ in \secref{sec:fabricius}
is a tridiagonal matrix as a matrix with cyclic indices
\eqref{S:matrix:fab}.

In \secref{sec:fabricius} we made a slight and technical modification to
Fabricius's construction by using the results in our old work
\cite{tak:92}, \cite{tak:95}. This modification gives an interpretation
of Baxter's strategy in \cite{bax:72} from the viewpoint of the
generalised algebraic Bethe Ansatz by Takhtajan and Faddeev
\cite{takh-fad:79}: Baxter used a ``half'' of the Gauss decomposition of
the gauge transfromation matrix in the generalised algebraic Bethe
Ansatz. See \remref{rem:gauss-decomp-Mlambda}.

Sklyanin's unitary operators $U_1$ and $U_3$ commute with $Q$-operators
constructed in \secref{sec:baxter72} and \secref{sec:fabricius} up to
shifts of the spectral parameter and a constant multiplication, as is
proved in \secref{subsec:quasi-periodicity:baxter72} and
\secref{subsec:quasi-periodicity:fabricius}.

We make several concluding remarks with comments on related works in the
final section \secref{sec:conclusion}.

All the necessary facts about the Sklyanin algebra are collected in
\appref{app:sklyanin}. In the proof of the commutation relation of $Q_R$
and $Q_L$ we need to compute the Sklyanin form (a scalar product in the
space of theta functions) of several vectors, using Rosengren's results
\cite{ros:04} and \cite{ros:07}. We summarise necessary results in
\appref{app:sklyanin-form}.

\subsection*{Notations}

Throughout this paper we use the following notations and symbols.
\begin{itemize}
 \item $N\in \Integer_{>0}$: the number of sites.
 \item $l\in \frac{1}{2}\Integer_{>0}$: the spin of the representation at
       each site. We assume that $Nl\in\Integer$. Namely, if $l$ is a
       half integer, then $N$ should be even.
 \item $\tau \in i \Real_{>0}$; the elliptic modulus, which is purely
       imaginary. 
 \item $\eta=\dfrac{r'}{2lr} \in
       \left[-\dfrac{1}{2(2l+1)},\dfrac{1}{2(2l+1)}\right]$
       ($r,r'\in\Integer$, $r>0$): the anisotropy parameter.
 \item The notations for the theta functions are the same as those in
       Mumford's book \cite{mum:83} (cf.\ Sklyanin's papers
       \cite{skl:82}, \cite{skl:83}):
\begin{equation}
    \theta_{ab}(z,\tau) = \sum_{n\in\Integer}
       \exp\left( 
             \pi i \left( \frac{a}{2} + n \right)^2 \tau
           +2\pi i \left( \frac{a}{2} + n \right)
                   \left( \frac{b}{2} + z \right)
           \right).
\label{def:theta}
\end{equation}
       (cf.\ Jacobi's notation (e.g., \cite{whi-wat}): 
       $\vartheta_1(\pi z,\tau)=-\theta_{11}(z,\tau)$, 
       $\vartheta_2(\pi z,\tau)=\theta_{10}(z,\tau)$, 
       $\vartheta_3(\pi z,\tau)=\theta_{00}(z,\tau)$, 
       $\vartheta_4(\pi z,\tau)=\theta_{01}(z,\tau)$.)
 \item We denote $\theta_{11}(z,\tau)$ by $[z]$ for simplicity. 
 \item $[z]_k := \prod_{j=0}^{k-1}[z+2j\eta] =
       [z][z+2\eta]\dotsb[z+2(k-1)\eta]$ for $k=1,2,\dotsc$,
       $[z]_0=1$. The following function is important.
\begin{equation}
    [z;a]_k := [z+a]_k [-z+a]_k.
\label{def:[z;a]}
\end{equation}
 \item The Pauli matrices are defined as usual:
\begin{equation*}
    \sigma^0 = \begin{pmatrix} 1 & 0  \\  0  & 1 \end{pmatrix},\quad
    \sigma^1 = \begin{pmatrix} 0 & 1  \\  1  & 0 \end{pmatrix},\quad
    \sigma^2 = \begin{pmatrix} 0 & -i \\  i  & 0 \end{pmatrix},\quad
    \sigma^3 = \begin{pmatrix} 1 & 0  \\  0  &-1 \end{pmatrix}.
\end{equation*}
\end{itemize}

\section{Definition of the model and the $Q$-operator}
\label{sec:model}

In this section we review the higher spin generalisation of the eight
vertex model, which was studied in \cite{tak:92}, \cite{tak:95},
\cite{tak:16}. Normalisations are the same as in \cite{tak:16}. For
notations see \appref{app:sklyanin}.

The local state spaces are the spin $l$ representation space of the
Sklyanin algebra ($l\in\frac{1}{2}\Integer_{>0}$) $V_i$
($i=1,\dotsc,N$): $V_i \cong \Theta_{00}^{4l+}$.  The total Hilbert
space $\calH$ is their tensor product:
\begin{equation}
    \calH := V_N \tensor \dotsb \tensor V_1,
\label{def:H}
\end{equation}
while the auxiliary space $V_0$ is a two-dimensional space:
$V_0\cong\Comp^2$. 

The transfer matrix $T(u)$ of the model is defined by
\begin{equation}
    T(u) :=
    \tr_0 L_N(u) L_{N-1}(u) \dotsb L_1(u),
\label{def:T(u)}
\end{equation}
where $\tr_0$ is the trace over the space $V_0$ and the $L$-operator
$L_n(u) \in \End_\Comp(\calH\tensor V_0)$ is defined by a representation
of the Sklyanin algebra on $\calH$,
$
    \rho_n := 1 \tensor \dotsb \tensor
    \overset{n\text{-th component}}{\rho^l} \tensor \dotsb \tensor 1,
$
as
\begin{equation}
    L_n(u) = \sum_{a=0}^3 W^L_a(u) \rho_n(S^a) \tensor \sigma^a.
\label{def:Ln}
\end{equation}

When the spin $l$ is $1/2$, the transfer matrix $T(u)$ is essentially
that of the eight vertex model. As usual, the $RLL$-relation \eqref{RLL}
leads to the commutativity of the transfer matrix:
\begin{equation}
    T(u) T(u') = T(u') T(u).
\label{TT}
\end{equation}

As in the previous paper \cite{tak:16}, our goal in this article is to
construct a {\em $Q$-operator} $Q(u)$, an operator satisfying the
commutation relations:
\begin{align}
    T(u) Q(u) &= h_-(u) Q(u-2\eta) + h_+(u) Q(u+2\eta),
\label{TQ}
\\
    Q(u) T(u) &= h_-(u) Q(u-2\eta) + h_+(u) Q(u+2\eta),
\label{QT}
\\
   Q(u) Q(v) &= Q(v) Q(u),
\label{QQ=QQ}
\end{align}
where the functions $h_\pm(u)$ are defined by
\begin{equation}
    h_\pm(u) := (2[u\mp2l\eta])^N.
\label{def:h+-}
\end{equation}
We apply two methods to this model, one by Baxter \cite{bax:72} and the
other by Fabricius \cite{fab:07} (see also \cite{roa:07-2}). Their main
strategies are essentially the same, which is due to Baxter
\cite{bax:72}. 
\begin{enumerate}
 \item First we construct the ``$Q_R$-operator'' of the form
\begin{equation}
    Q_R(u)
    = \tr_{\Comp^r}S_N(u) S_{N-1}(u) \dotsb S_1(u),
\label{QR=tr(SSS)}
\end{equation}
       satisfying the $TQ$-relation \eqref{TQ}. The operator $S_n(u)\in
       \Hom_{\Comp}((\Comp^{2l+1})^{\tensor N}\tensor \Comp^{r_0},
       \calH\tensor \Comp^{r_0})$ is an $r_0\times r_0 $-matrix of
       linear maps $(\Comp^{2l+1})^{\tensor N}\to \calH$:
\begin{equation}
    S_n(u) = \bigl((S^{i''}_{j''}(u))_n\bigr)_{i'',j''=1,\dotsc,r_0},
\label{Sn(u)}
\end{equation}
       where 
\begin{equation*}
    (S^{i''}_{j''}(u))_n
    =
    \Id_{V_N} \tensor \dotsb \tensor
    \overset{n\text{-th component}}{S^{i''}_{j''}(u)} \tensor
    \dotsb \tensor \Id_{V_1},
\end{equation*}
       and $S^{i''}_{j''}(u): \Comp^{2l+1}\to \Theta^{4l+}_{00}$. We
       denote
\begin{equation}
    S(u) = \bigl(S^{i''}_{j''}(u)\bigr)_{i'',j''=1,\dotsc,r_0}.
\label{S(u)}
\end{equation}
       The matrix size $r_0$ for the first method is $r$, while $r_0=2r$
       for the second method. (Recall that $r$ is a factor of the
       denominator of $\eta=\dfrac{r'}{2lr}$.) 
      
 \item The product $T(u) Q_R(u)$ has the form
\begin{equation*}
    T(u) Q_R(u)
    =
    \tr_{V_0\tensor \Comp^r} \prod_{1\leq n \leq N}^{\curvearrowleft}
    (L_n(u)\tensor S_n(u)),
\end{equation*}
       where the tensor product of $L_n(u)$ and $S_n(u)$ is a
       $(2r_0\times 2r_0)$-matrix with
       $\Hom_\Comp((\Comp^{2l+1})^{\tensor N},\calH)$-valued elements,
       $L^{i'}_{j'}(u)\, S^{i''}_{j''}(u)$.  We index the components of
       $V_0\cong\Comp^2$ by $i',j'=\pm$ and $\Comp^{r_0}$ by $i'', j''
       \in \{1,\dotsc,r_0\}$.

       If there exists a $(2r_0\times 2r_0)$-matrix $M$ with
       scalar elements, such that
\begin{equation}
    M^{-1} (L(u)\tensor S(u)) M
    =
    \begin{pmatrix}
     A(u) & 0 \\ * & D(u)
    \end{pmatrix} \text{ or }
    \begin{pmatrix}
     A(u) & * \\ 0 & D(u)
    \end{pmatrix},
\label{M-1(LS)M}
\end{equation}
       then we have
\begin{equation}
    T(u) Q_R(u)
    = \tr_{\Comp^{r_0}} A_N(u) \dotsb A_1(u)
    + \tr_{\Comp^{r_0}} D_N(u) \dotsb D_1(u).
\label{TQR:temp}
\end{equation}
       We derive the $TQ$-relation \eqref{TQ} from \eqref{TQR:temp} by
       using explicit forms of $A(u)$ and $D(u)$.
 \item The $Q_L$-operator which satisfies the $QT$-relation \eqref{QT}
       is obtained as an adjoint to $Q_R$ with a change of the spectral
       parameter. 
 \item The commutativity $Q_L(u)Q_R(u')=Q_L(u')Q_R(u)$ is shown by
       computing the matrix elements of its both sides explicitly.
 \item The $Q$-operator is defined by
       $Q(u)=Q_R(u)Q_R(u_0)^{-1}=Q_L(u_0)^{-1}Q_L(u)$, where $u_0$ is
       chosen so that $Q_R(u_0)$ and $Q_L(u_0)$ are invertible. The
       existence of such $u_0$ will not be discussed in this paper. (For
       the eight vertex model, detailed study of degeneracy of $Q_R(u)$
       is in \cite{fab-mcc:03-05}.)

       The commutation relations \eqref{TQ}, \eqref{QT} and
       \eqref{QQ=QQ} follow automatically from \eqref{TQ} for $Q_R$,
       \eqref{QT} for $Q_L(u)$ and $Q_L(u)Q_R(u')=Q_L(u')Q_R(u)$.
 \item The $Q$-operator is quasi-periodic in the spectral parameter
       $u$. The shift of $u$ is equivalent to multipication with
       involutions $U_1^{\tensor N}$ and $U_3^{\tensor N}$ on the
       Hilbert space $\calH$. 
\end{enumerate}

\section{The $Q$-operator {\`a} la Baxter (1972)}
\label{sec:baxter72}

In this section we construct the $Q$-operator for the higher spin eight
vertex model, following the strategy sketched in \secref{sec:model} with
$r_0=r$. In this construction, the matrix $S(u)$ \eqref{S(u)} is of the
form,
\begin{equation}
    S(u)=
     \begin{pmatrix}
     S^1_1(u) & S^1_2(u) &          &              &
     \\
     S^2_1(u) & 0        & S^2_3(u) &              &
     \\
              & S^3_2(u) & 0        & \ddots       &
     \\
              &          & \ddots   &              & S^{r-1}_r(u)
     \\
              &          &          & S^r_{r-1}(u) & S^r_r(u)
    \end{pmatrix}.
\label{S:matrix:bax}
\end{equation}

\subsection{Construction of $Q_R$}
\label{subsec:QR:baxter72}

As in Appendix C of \cite{bax:72}, we may assume that the matrix $M$ in
\eqref{M-1(LS)M} has the form
\begin{equation}
    M = \sum_{i''=1}^r \tilde M_{i''} \tensor E^{i''}_{i''},
\label{M:general}
\end{equation}
where
$E^{i''}_{j''}=(\delta^{i'',k''}\delta_{j'',l''})_{k'',l''=1,\dotsc,r}$
is the $r\times r$-matrix unit and $\tilde M_{i''}$ is an upper
triangular $2\times 2$-matrix,
\begin{equation}
    \tilde M_{i''}
    =
    \begin{pmatrix}
     1 & p_{i''} \\ 0 & 1
    \end{pmatrix}.
\label{tildeM}
\end{equation}
Note that the inverse of $M$ is of the following form:
\begin{equation}
    M^{-1} = \sum_{i''=1}^r \tilde M_{i''}^{-1} \tensor E^{i''}_{i''}.
\label{M:general}
\end{equation}
The problem of constructing $S$ and $M$ satisfying \eqref{M-1(LS)M}
reduces to finding
\begin{itemize}
 \item the element $p_{i''}$ of $\tilde M_{i''}$ \eqref{tildeM} such
       that the (1,2)-element of $\tilde M_{i''}^{-1}L(u)\tilde M_{j''}$
       is degenerate.
 \item an operator $S^{i''}_{j''}$, the image of which belongs to the
       kernel of the (1,2)-element of $\tilde M_{i''}^{-1}L(u)\tilde
       M_{j''}$.
\end{itemize}

Let us denote the elements of the $L$-operator as follows:
\begin{equation}
    L(u) = \sum_{a=0}^3 W^L_a(u) \rho^l(S^a) \tensor \sigma^a
    =
    \begin{pmatrix}
     L^-_-(u) & L^-_+(u) \\ L^+_-(u) & L^+_+(u)
    \end{pmatrix}.
\label{L}
\end{equation}
Then, the elements of
\begin{equation}
    \tilde M_{i''}^{-1} L(u) \tilde M_{j''}
    =
    \begin{pmatrix}
     \alpha_{i'',j''}(u) & \beta_{i'',j''}(u)\\
     \gamma_{i'',j''}(u) & \delta_{i'',j''}(u)
    \end{pmatrix}
\label{MLM-1}
\end{equation}
have the following form:
\comment{cf. note22 p.20 (1)20/II/2017}
\begin{align}
    \alpha_{i'',j''} &=
    L^-_-(u) - p_{i''}L^+_-,
\label{alpha-ij}
\\
    \beta_{i'',j''} &=
    L^-_-(u) p_{j''} + L^-_+ - p_{i''}L^+_- p_{j''} - p_{i''}L^+_+,
\label{beta-ij}
\\
    \gamma_{i'',j''} &=
    L^+_-(u),
\label{gamma-ij}
\\
    \delta_{i'',j''} &=
    L^+_-(u) p_{j''} + L^+_+.
\label{delta-ij}
\end{align}
Substituting the expression \eqref{L}, we have
\comment{note22 p21, (2)20/II/2017}
\begin{equation}
 \begin{split}
    \beta_{i'',j''} ={}&
     (-p_{i''}+p_{j''}) W_0(u)\rho^{l}(S_0)
    +(1-p_{i''}p_{j''}) W_1(u)\rho^{l}(S_1)
\\
    &+(-1-p_{i''}p_{j''})W_2(u)\rho^{l}(S_2)
     +(p_{i''}+p_{j''})  W_3(u)\rho^{l}(S_3).
 \end{split}
\label{beta-ij:S}
\end{equation}
When we choose $p_{i''}$ and $p_{j''}$ as
\comment{note23 p.10,(2)4/VII/2017}
\begin{equation}
    p_{i''} 
    =
     -\frac
      {\theta_{00}\bigl(\frac{\lambda\pm4l\eta}{2},\frac{\tau}{2}\bigr)}
      {\theta_{01}\bigl(\frac{\lambda\pm4l\eta}{2},\frac{\tau}{2}\bigr)},
    \qquad
    p_{j''} 
    =
     -\frac
      {\theta_{00}\bigl(\frac{\lambda}{2},\frac{\tau}{2}\bigr)}
      {\theta_{01}\bigl(\frac{\lambda}{2},\frac{\tau}{2}\bigr)}
\label{p(lambda)}
\end{equation}
($\lambda\in\Comp$ is a parameter), the operator $\beta_{i'',j''}$
\eqref{beta-ij:S} is degenerate with a null vector
\comment{note23 p.10 (2)4/VII/2017 (with an error: lambda should be
lambda/2)}
\comment{(cf. p.3 (1)29/VI/2017, p.30 (2)27/X/2017) }
\comment{note25 p.4 (1)11/III/2018}
\begin{equation}
    f_\pm(\lambda,u,z)
    :=
    \left[z;\frac{\pm\lambda+u}{2}+(-l+1)\eta \right]_{2l}
    \in\Theta^{4l++}_{00},
\label{null-gamma(lambda)}
\end{equation}
where the notation $[z;a]_{2l}$ is defined by \eqref{def:[z;a]}. This
vector is transformed by the actions of $\alpha_{i'',j''}$ and
$\delta_{i'',j''}$ as follows.
\begin{align}
\comment{note24 p.13 (3)2/I/2018}
    \alpha_{i'',j''}f_\pm(\lambda,u,z)
    &=
     \frac{\theta_{01}\bigl(\frac{\lambda}{2},\frac{\tau}{2}\bigr)}
    {\theta_{01}\bigl(\frac{\lambda\pm4l\eta}{2},\frac{\tau}{2}\bigr)}
    2\theta_{11}(u-2l\eta)
    f_\pm(\lambda,u+2\eta,z),
\label{alpha-f}
\\
\comment{note24 p.21 (6)3/I/2018}
    \delta_{i'',j''}f_\pm(\lambda,u,z)
    &=
    \frac
    {\theta_{01}\bigl(\frac{\lambda\pm4l\eta}{2},\frac{\tau}{2}\bigr)}
    {\theta_{01}\bigl(\frac{\lambda}{2},\frac{\tau}{2}\bigr)}
    2\theta_{11}(u+2l\eta)
    f_\pm(\lambda,u-2\eta,z).
\label{delta-f}
\end{align}

Let us set $\lambda=2(2j''-1)l\eta$:
\comment{note24 p.22 (1) 26/I/2018}
\begin{equation}
    p_{j''} 
    :=
     -\frac
      {\theta_{00}\bigl((2j''-1)l\eta,\frac{\tau}{2}\bigr)}
      {\theta_{01}\bigl((2j''-1)l\eta,\frac{\tau}{2}\bigr)}.
\label{pi''}
\end{equation}
Then, as we mentioned above, $\beta_{j''\pm1,j''}$ is degenerate. In
addition $\beta_{1,1}$ and $\beta_{r,r}$ are also degenerate, since
$p_0=p_1$ and $p_r=p_{r+1}$ because of the evenness and the periodicity
of theta functions. (Recall $2lr\eta=r'\in\Integer$.) It is easy to
check that \eqref{alpha-f} and \eqref{delta-f} hold not only for
$i''=j''\pm1$ but also for $i''=j''=1$ and $i''=j''=r$.

Therefore, defining an operator
$S^{i''}_{j''}(u):\Comp^r\to\Theta^{4l++}_{00}$ by
\comment{note24 p.24 (2)27/I/2018, note25 p.1 (1)27/II/2018}
\begin{equation}
 \begin{aligned}
    S^{j''\pm1}_{j''}(u)&: \Comp^r\owns e_k \mapsto
    \tau_{kj''} f_\pm(2(2j''-1)l\eta,u,z)&& \in \Theta^{4l++}_{00},
\\
    S^{1}_{1}(u)&: \Comp^r\owns e_k \mapsto
    \tau_{k1} f_-(2l\eta,u,z)&& \in \Theta^{4l++}_{00},
\\
    S^{r}_{r}(u)&: \Comp^r\owns e_k \mapsto
    \tau_{kr} f_+(2(2r-1)l\eta,u,z)&& \in \Theta^{4l++}_{00},
 \end{aligned}
\label{Si''j'':bax}
\end{equation}
($\{e_k\}_{k=1,\dotsc,r}$ is the standard basis of $\Comp^r$) with
generic parameters $\tau_{kj''}$ and
\begin{equation}
    S^{i''}_{j''}(u) = 0
    \text{ unless } i''=j''\pm1,\ i''=j''=1 \text{ or } i''=j''=r,
\label{Sn(u):zero-element:bax}
\end{equation}
(cf.\ \eqref{S:matrix:bax}), we obtain the operator-valued matrix $S(u)$
\eqref{S(u)} satisfying \eqref{M-1(LS)M}.

It follows from \eqref{alpha-f} and \eqref{delta-f} that the operators
$A(u)$ and $D(u)$ in \eqref{M-1(LS)M} are
\begin{equation}
 \begin{aligned}
    A(u)
    &=
    X^{-1} \bigl(2[u-2l\eta] S(u+2\eta)\bigr) X,
\\
    D(u)
    &=
    X      \bigl(2[u+2l\eta] S(u-2\eta)\bigr) X^{-1},
 \end{aligned}
\label{A,D:bax}
\end{equation}
where $X$ is a diagonal matrix
$
    X
    =
    \diag_{i''=1\dotsc,r}\Bigl(
    \theta_{01}\bigl((2i''-1)l\eta,\frac{\tau}{2}\bigr)
    \Bigr).
$

Thus the equation \eqref{TQR:temp} reduces to the $TQ$-relation,
\begin{equation}
     T(u) Q_R(u) = h_-(u) Q_R(u-2\eta) + h_+(u) Q_R(u+2\eta)
\label{TQR}
\end{equation}
with the functions $h_\pm(u)$ defined by \eqref{def:h+-}.

\subsection{Hermitian conjugate and $Q_L$}
\label{subsec:QL:baxter72}

The operator $Q_L(u)$  satisfying the  ``$QT$''-relation \eqref{QT}  is
constructed in a similar way as in \cite{tak:16} with the help of the
Sklyanin form \eqref{skl-form}. 

The proof of the equation (3.24) in \cite{tak:16} shows that the adjoint
operator of the transfer matrix $T(u)$ with respect to the Sklyanin form
on $\calH$ satisfies
\begin{equation}
 \begin{split}
    \bigl(T(u)\bigr)^* = (-1)^N T(-\bar u).
\end{split}
\label{T*}
\end{equation}
(When $N$ is even, this is (3.24) in \cite{tak:16}.) Therefore, if we
define the Hermitian structure of $\Comp^r$ by $(e_i,e_j)=\delta_{ij}$,
the adjoint of \eqref{TQR} with $-\bar u$ instead of $u$ gives
\begin{equation*}
    Q_R(-\bar u)^* T(u) 
    =
    h_+(u) Q_R(-\bar u - 2\eta)^*
    +
    h_-(u) Q_R(-\bar u + 2\eta)^*.
\end{equation*}
(Note that $h_{\pm}(-\bar u)=(-1)^N \overline{h_{\pm}(u)}$.) Thus the
operator $Q_L(u)$ defined by
\begin{equation}
    Q_L(u) := Q_R(-\bar u)^*:\calH \to (\Comp^{2l+1})^{\tensor N}
\label{def:QL}
\end{equation}
satisfies the $QT$-relation \eqref{QT}:
\begin{equation}
    Q_L(u) T(u) = h_-(u) Q_L(u-2\eta) + h_+(u) Q_L(u+2\eta).
\label{QLT}
\end{equation}

\subsection{Commutation relation of $Q_R$ and $Q_L$}
\label{subsec:QLQR=QLQR:baxter72}

Next we prove the commutation relation of $Q_R(u)$ and $Q_L(u)$,
\begin{equation}
    Q_L(u) Q_R(u') = Q_L(u') Q_R(u).
\label{QLQR=QLQR}
\end{equation}
Since the operator $Q_R(u)$ acts on the basis vector of $\Comp^r$ as
\[
    Q_R(u) e_{i_N}\tensor \dotsb \tensor e_{i_1}
    =
    \sum_{i''_1,\dotsc,i''_N\in\{1,\dots,r\}}
    S^{i''_N}_{i''_{N-1}}(u) e_{i_N} \tensor \dotsb \tensor
    S^{i''_1}_{i''_{N}}(u)   e_{i_1},
\]
the matrix elements of $Q_L(u) Q_R(u')$ has the form
\[
 \begin{split}
    &(              e_{j_N}\tensor \dotsb \tensor e_{j_1},
     Q_L(u) Q_R(u') e_{i_N}\tensor \dotsb \tensor e_{i_1})
\\
    ={}
    &\langle Q_R(-\bar u) e_{j_N}\tensor \dotsb \tensor e_{j_1},
             Q_R(u')      e_{i_N}\tensor \dotsb \tensor e_{i_1}\rangle
\\
    ={}
    &\sum_{\substack{i''_1,\dotsc,i''_N, j''_1,\dotsc,j''_N\\ \in
           \{1,\dotsc,r\}}}
    \prod_{k=1}^N
    \langle S^{j''_k}_{j''_{k-1}}(-\bar u) e_{j_k},
            S^{i''_k}_{i''_{k-1}}(     u') e_{i_k} \rangle
\\
    ={}&
    \tr_{\Comp^r\tensor \Comp^r}
    W(j_N,i_N| u,u') \dotsb W(j_1,i_1| u,u'),
 \end{split}
\]
where the $r^2\times r^2$-matrix $W(j,i|u,v)$ is defined by
\begin{equation}
    W(j,i|u,v)
    =
    \left(
     \langle S^{j'''}_{j''}(u) e_{j},
             S^{i'''}_{i''}(v) e_{i} \rangle
    \right)_{(j''',i'''), (j'',i'')\in \{1,\dotsc,r\}^2}.
\label{def:W(ji|uv)}
\end{equation}
Hence, as in Appendix C of \cite{bax:72} (see also \cite{fab-mcc:07}),
if we can find a matrix
$Y=(Y^{j'',i''}_{j''',i'''})_{(j'',i''),(j''',i''')\in\{1,\dotsc,r\}^2}$
such that 
\begin{equation}
    Y W(j,i|u,v) Y^{-1} = W(j,i|v,u),
\label{YWY-1}
\end{equation}
all the matrix elements of $Q_L(u) Q_R(u')$ are symmetric in $(u,u')$,
which proves the commutation relation \eqref{QLQR=QLQR}.

Using the results by Rosengren, \cite{ros:04} and \cite{ros:07} (see
\appref{app:sklyanin-form}), we can compute the non-zero elements of
$W(j,i|u,v)$ explicitly:
\begin{multline}
    W(j,i|u,v)^{j''',i'''}_{j'',i''}
    =\langle S^{j'''}_{j''}(u) e_{j},S^{i'''}_{i''}(v) e_{i} \rangle
\\
    =
    \overline{\tau_{jj''}}\tau_{ii''}
    F\left(\frac{u+v}{2}+ 2 {w^{(F)}}^{j'''i'''}_{j''i''}l\eta \right)
    G\left(\frac{v-u}{2}+ 2 {w^{(G)}}^{j'''i'''}_{j''i''}l\eta \right),
\label{W(j,i|u,v):element:bax}
\end{multline}
The functions $F$ and $G$ are defined by \eqref{def:F} and \eqref{def:G}
in \appref{app:sklyanin-form} respectively and constants
${w^{(F)}}^{j''i''}_{j'''i'''}$ and ${w^{(G)}}^{j''i''}_{j'''i'''}$ are
defined in \tableref{tab:wFjiji} and \tableref{tab:wGjiji}
respectively. Except for sixteen cases in these tables,
$W(j,i|u,v)^{j'',i''}_{j''',i'''}$ vanishes because of the condition
\eqref{Sn(u):zero-element:bax}.

\begin{table}[h]
\comment{note25 pp.2-24 5-27/III/2018}
\begin{center}
 \begin{tabular}{|l||c|c|c|c|}
 \hline
  ${w^{(F)}}^{j'''i'''}_{j''i''}$
  & $j'''=j''+1$ & $j'''=j''-1$ & $j'''=j''=1$ & $j'''=j''=r$ \\
 \hline\hline
  $i'''=i''+1$ & $i''-j''$ & $i''+j''-1$ & $i''$ & $i''-r$ \\
 \hline
  $i'''=i''-1$ & $-i''-j''+1$ & $-i''+j''$ & $-i''+1$ & $-i''-r+1$  \\
 \hline
  $i'''=i''=1$ & $-j''$ & $j''-1$ & $0$ & $r$ \\
 \hline
  $i'''=i''=r$ & $-j''+r$ & $j''+r-1$ & $-r$ & $0$ \\
 \hline
 \end{tabular}
\caption{${w^{(F)}}^{j'''i'''}_{j''i''}$ in
 \eqref{W(j,i|u,v):element:bax}.} 
 \label{tab:wFjiji}
\end{center}
\end{table}

\begin{table}[h]
\comment{note25 p.13 (3) 14/III/2018}
\comment{argument of G shifted by $2(-2l+1)l\eta$}
\begin{center}
 \begin{tabular}{|l||c|c|c|c|}
 \hline
  ${w^{(G)}}^{j'''i'''}_{j''i''}$
  & $j'''=j''+1$ & $j'''=j''-1$ & $j'''=j''=1$ & $j'''=j''=r$ \\
 \hline\hline
  $i'''=i''+1$ & $i''+j''$ & $i''-j''+1$ & $i''$ & $i''+r$ \\
 \hline
  $i'''=i''-1$ & $-i''+j''+1$ & $-i''-j''+2$ & $-i''+1$ & $-i''+r+1$ \\
 \hline
  $i'''=i''=1$ & $j''$ & $-j''+1$ & $0$ & $r$ \\
 \hline
  $i'''=i''=r$ & $j''+r$ & $-j''+r+1$ & $r$ & $2r$ \\
 \hline
 \end{tabular}
\caption{${w^{(G)}}^{j'''i'''}_{j''i''}$ in
 \eqref{W(j,i|u,v):element:bax}.} 
 \label{tab:wGjiji}
\end{center}
\end{table}

As in Appendix C of \cite{bax:72}, we take a diagonal matrix $Y$ in
\eqref{YWY-1} of the following form:
\comment{note25 p.26 (1) 29/III/2018}
\begin{equation}
    Y= \diag_{(j''i'')\in\{1,\dotsc,r\}^2} (y_{j''i''}),
    \quad
    y_{j''i''} = t_{j''+i''} t_{-j''+i''+1},
\label{def:Y}
\end{equation}
where $t_m$ is defined by the recurrence relation:
\begin{equation}
    \frac{t_{m+2}}{t_m}
    =
    \frac{G(\frac{u-v}{2}+2ml\eta)}{G(\frac{v-u}{2}+2ml\eta)}.
\label{tm:recurrence}
\end{equation}
We have to show \eqref{YWY-1}, which is equivalent to
\begin{equation}
    \frac{y_{j'''i'''}}{y_{j''i''}}
    =
    \frac{G(\frac{u-v}{2}+ 2 {w^{(G)}}^{j'''i'''}_{j''i''}l\eta )}
         {G(\frac{v-u}{2}+ 2 {w^{(G)}}^{j'''i'''}_{j''i''}l\eta )}
\label{y/y=G/G}
\end{equation}
by virtue of \eqref{W(j,i|u,v):element:bax} and
\eqref{def:Y}. Non-trivial cases of the left hand side of
\eqref{y/y=G/G} expressed in $t_m$'s are listed in \tableref{tab:y/y}.

\begin{table}[h]
\comment{note25 p.26 (1) 29/III/2018}
\comment{$y_{i''j''}$ in the note is $y_{j''i''}$ here.}
\begin{center}
 \begin{tabular}{|l||c|c|c|c|}
 \hline
  $y_{j'''i'''}/y_{j''i''}$
  & $j'''=j''+1$ & $j'''=j''-1$ & $j'''=j''=1$ & $j'''=j''=r$ \\
 \hline\hline
  $i'''=i''+1$ & $\frac{t_{i''+j''+2}}{t_{i''+j''}}$
               & $\frac{t_{i''-j''+3}}{t_{i''-j''+1}}$
	       & $\frac{t_{i''+2}}{t_{i''}}$
               & $\frac{t_{i''+r+1}\,t_{i''-r+2}}{t_{i''+r}\,t_{i''-r+1}}$ \\
 \hline
  $i'''=i''-1$ & $\frac{t_{i''-j''-1}}{t_{i''-j''+1}}$
               & $\frac{t_{i''+j''-2}}{t_{i''+j''}}$
	       & $\frac{t_{i''-1}}{t_{i''+1}}$
               & $\frac{t_{i''+r-1}\,t_{i''-r}}{t_{i''+r}\,t_{i''-r+1}}$ \\
 \hline
  $i'''=i''=1$ & $\frac{t_{j''+2}\,t_{-j''+1}}{t_{j''+1}\,t_{-j''+2}}$
               & $\frac{t_{j''}\,t_{-j''+3}}{t_{j''+1}\,t_{-j''+2}}$
	       & 1
               & 1 \\
 \hline
  $i'''=i''=r$ & $\frac{t_{j''+r+1}\,t_{-j''+r}}{t_{j''+r}\,t_{-j''+r+1}}$
               & $\frac{t_{j''+r-1}\,t_{-j''+r+2}}{t_{j''+r}\,t_{-j''+r+1}}$
	       & 1
               & 1 \\
 \hline
 \end{tabular}
\caption{$y_{j'''i'''}/y_{j''i''}$.}
 \label{tab:y/y}
\end{center}
\end{table}

Comparing \tableref{tab:wGjiji} and \tableref{tab:y/y} and using
properties of the function $G$, \eqref{G:even}\footnote{From the
evenness \eqref{G:even} of $G$ follows $t_{m+1}=t_{-m+1}$ by
induction. This property is useful when one checks \eqref{y/y=G/G}.} and
\eqref{G:periodic}, we can show \eqref{y/y=G/G}. Thus the commutation
relation \eqref{QLQR=QLQR} has been proved.

\medskip
If there exists $u_0\in\Comp$ such that $Q_R(u_0)$ and $Q_L(u_0)$ are
non-degenerate, multiplying $Q_R(u_0)^{-1}$ from the right and
$Q_L(u_0)^{-1}$ from the left to \eqref{QLQR=QLQR} ($u'=u_0$), we obtain
the $Q$-operator
\begin{equation}
    Q(u):= Q_R(u) Q_R(u_0)^{-1} = Q_L(u_0)^{-1} Q_L(u),
\label{Q=QRQR-1=QL-1QL}
\end{equation}
which satisfies the $TQ$- and $QT$-relations, \eqref{TQ} and \eqref{QT},
because of \eqref{TQR} and \eqref{QLT}. Commutativity \eqref{QQ=QQ} is a
direct consequence of \eqref{QLQR=QLQR}.

Non-degeneracy of $Q_R(u)$ and $Q_L(u)$ is quite non-trivial, as is
discussed in \cite{fab-mcc:03-05} (see also \cite{fab-mcc:07},
\cite{fab-mcc:09}), which showed that $Q_R(u)$ has a large kernel for
certain values of $\eta$ and $N$. This problem would be discussed in a
separate paper.

\subsection{Quasi-periodicity of the $Q$-operator}
\label{subsec:quasi-periodicity:baxter72}

The unitary operators $U_1$ and $U_3$, \eqref{def:Ua}, act on
$f_\pm(\lambda,u,z)$, \eqref{null-gamma(lambda)}, as follows.
\begin{equation}
 \begin{aligned}
\comment{note26 p.2 (1)23/V/2018}
    U_1 f_\pm(\lambda,u,z) &= e^{-l \pi i} f_\pm(\lambda,u+1,z),
\\
\comment{note26 p.4 (2)17/VI/2018}
    U_3 f_\pm(\lambda,u,z)
    &= e^{l(\tau-1)\pi i +2l(\pm\lambda+u+2l\eta)\pi i }
       f_\pm(\lambda,u+\tau,z).
 \end{aligned}
\label{Uaf+-}
\end{equation}
Hence, acting $U_a$ on \eqref{Si''j'':bax}, we obtain
\comment{note26 p.9 (2)26/VI/2018}
\begin{equation}
 \begin{aligned}
    U_1 S^{i''}_{j''}(u)&= e^{-l\pi i} S^{i''}_{j''}(u+1),
\\
    U_3 S^{i''}_{j''}(u)&=
    e^{l(\tau-1)\pi i +2l \pi i u+ 8\, d(i'',j'')\, l^2 \eta\pi i}
    S^{i''}_{j''}(u+\tau),
 \end{aligned}
\label{UaS(u):bax}
\end{equation}
where
\begin{equation}
    d(i'',j'')=
    \begin{cases} j'' &(i''=j''+1), \\ -i'' &(i''=j''-1), 
      \\ 0 &(i''=j''=1 \text{ or }r).
    \end{cases}
\label{def:d(i'',j''):bax}
\end{equation}
(Recall that $8l^2r\eta\in2\Integer$.) The quasi-periodicity of $Q_R(u)$
with respect to $u\mapsto u+1$ is derived directly from \eqref{UaS(u):bax},
because of the definitions \eqref{QR=tr(SSS)} and \eqref{Sn(u)}:
\comment{note26 p.6 (1)24/VI/2018}
\begin{equation}
    U_1^{\tensor N} Q_R(u) = e^{-N\pi i l} Q_R(u+1).
\label{U1QR(u)=QR(u+1):bax}
\end{equation}
In order to derive the quasi-periodicity with respect to $u\mapsto
u+\tau$, a similarity transformation of $S(u)$ by a diagonal matrix 
\begin{equation}
    S(u) \mapsto A S(u) A^{-1},\quad
    A=\diag_{i''=1,\dotsc,r}(A_{i''}),\
    A_{i''} := \exp(4j''(j''-1)l^2\eta\pi i),
\label{def:A:bax}
\end{equation}
is necessary. This transformation does not change $Q_R(u)$ thanks to
$\tr_{\Comp^r}$ in \eqref{QR=tr(SSS)} and cancels the factors in
\eqref{UaS(u):bax} depending on $i''$ and $j''$. Thus we obtain 
\comment{note26 p.9 (2)26/VI/2018}
\begin{equation}
    U_3^{\tensor N} Q_R(u)
    =
    e^{Nl \pi i (\tau-1)+2 Nl\pi i u} Q_R(u+\tau).
\label{U3QR(u)=QR(u+tau):bax}
\end{equation}

The quasi-periodicity of $Q_L(u)$ is proved by taking the
Hermitian adjoint of relations \eqref{U1QR(u)=QR(u+1):bax} and
\eqref{U3QR(u)=QR(u+tau):bax}. For any vectors $\Phi\in\calH$ and
$v\in(\Comp^{2l+1})^{\tensor N}$, we have
\begin{equation}
    (Q_L(u) U_a^{\tensor N}\Phi,v)_{(\Comp^{2l+1})^{\tensor N}}
    =
    \langle
    \Phi, U_a^{\tensor N} Q_R(-\bar u)v
    \rangle_\calH.
\label{QLUa=UaQR} 
\end{equation}
by the definition \eqref{def:QL} of $Q_L(u)$ and the unitarity of $U_a$.

\comment{note27 p.28, (3)19/IX/2018}
As for the operator $U_1$, we have $U_1^{\tensor N} Q_R(-\bar
u)=e^{-N\pi i l} Q_R(-\bar u+1)$ because of \eqref{U1QR(u)=QR(u+1):bax}.
Note that $f_\pm(\lambda,u+2,z)=f_\pm(\lambda,u,z)$, which implies
$Q_R(u+2)=Q_R(u)$ by virtue of \eqref{QR=tr(SSS)} and
\eqref{Si''j'':bax}. Therefore the right hand side of
\eqref{QLUa=UaQR} is equal to
$
    \langle \Phi, e^{-N\pi i l} Q_R(-\bar u+1)v \rangle_\calH
    =
    \langle \Phi, e^{-N\pi i l} Q_R(-\bar u-1)v \rangle_\calH
$.
Taking the Hermitian adjoint of the operator again, we obtain $Q_L(u)
U_1^{\tensor N} = e^{N\pi il}Q_L(u+1)$. From the assumption
$Nl\in\Integer$ follows
\begin{equation}
    Q_L(u) U_1^{\tensor N} = e^{-Nl\pi i} Q_L(u+1).
\label{QL(u)U1=QL(u+1):bax}
\end{equation}

\comment{note27 p.27, (2)19/IX/2018}
As for the operator $U_3$, the composition $U_3^{\tensor N} Q_R(-\bar
u)$ in \eqref{QLUa=UaQR} becomes $e^{Nl\pi i(\tau-1)-2Nl\pi i\bar
u}Q_R(-\bar u+\tau)$. Because $\tau$ is purely imaginary, $-\bar
u+\tau=-\overline{(u+\tau)}$. The Hermitian adjoint is, by the
definition \eqref{def:QL},
$
    e^{Nl\pi i(\tau-1)+2Nl\pi i+2Nl\pi iu}Q_L(u+\tau).
$
Using the assumption $Nl\in\Integer$ again, we have
\begin{equation}
    Q_L(u) U_3^{\tensor N}
    =
    e^{Nl\pi i(\tau-1)+2Nl\pi iu} Q_L(u+\tau).
\label{QL(u)U3=QL(u+tau):bax}
\end{equation}
Since \eqref{U1QR(u)=QR(u+1):bax} and \eqref{QL(u)U1=QL(u+1):bax},
\eqref{U3QR(u)=QR(u+tau):bax} and \eqref{QL(u)U3=QL(u+tau):bax} have
pairwise the same coefficients in the right hand side, these
quasi-periodicity relations for $Q_R(u)$ and $Q_L(u)$ give the
quasi-periodicity of the $Q$-operator,
\begin{equation}
 \begin{aligned}
    U_1^{\tensor N} Q(u) &= 
    Q(u) U_1^{\tensor N}
    = e^{-Nl\pi i} Q(u+1),
\\
    U_3^{\tensor N} Q(u) &= 
    Q(u) U_3^{\tensor N}
    = e^{Nl\pi i(\tau-1)+2Nl\pi iu} Q(u+\tau).
 \end{aligned}
\label{Q:quasi-periodicity:bax}
\end{equation}

\section{The $Q$-operator {\`a} la Fabricius}
\label{sec:fabricius}

Here we construct the $Q$-operator of the same model in a different way,
following \cite{fab:07}. The strategy is almost the same as in
\secref{sec:baxter72}. The matrix $S(u)$ has the form
\begin{equation}
    S(u)=
     \begin{pmatrix}
     0           & S^1_2(u) &          &            & S^1_{2r}(u)
     \\
     S^2_1(u)    & 0        & S^2_3(u) &            &
     \\
                 & S^3_2(u) & 0        & \ddots     &
     \\
                 &          & \ddots   &            & S^{2r-1}_{2r}(u)
     \\
     S^{2r}_1(u) &          &          & S^{2r}_{2r-1}(u) & 0
    \end{pmatrix},
\label{S:matrix:fab}
\end{equation}
instead of \eqref{S:matrix:bax}\footnote{In \cite{fab:07} the size of
$S(u)$ is $r\times r$. In fact, for the construction of $Q_R$ satisfying
\eqref{TQR} and $Q_L$ satisfying \eqref{QLT} we can use $S(u)$ of size
$r\times r$ as in \cite{fab:07}. But when we prove the commutation
relation \eqref{QLQR=QLQR} in \secref{subsec:QL:fabricius}, we need
$S(u)$ of size $2r\times 2r$.}. Correspondingly, we use a different
matrix $M$, \eqref{M:general}.

Since the trace of a product of an odd number of matrices of the above
form is trivially zero, we have to assume that $N$ is even. 

\subsection{Construction of $Q_R$}
\label{subsec:QR:fabricius}

The essential idea in \secref{sec:baxter72} was to make one of
off-diagonal elements of $M^{-1}(L(u)\tensor S(u))M$ in \eqref{M-1(LS)M}
zero, which was a consequence of degeneracy of an off-diagonal element
of the twisted $L$-matrix $\tilde M^{-1}_{i''}L(u)\tilde M_{j''}$ in
\eqref{MLM-1}. We used such twisting (or gauge transformation) when we
applied Takhtajan-Faddeev's generalised Bethe Ansatz in
\cite{takh-fad:79} to the model (\cite{tak:92} and \cite{tak:95}), or
when we constructed the $Q$-operator by the method of Baxter's 1973
paper \cite{bax:73} (\cite{tak:16}).

In this section we use the matrix $M_\lambda(v)$ in (3.1) of
\cite{tak:16} (cf.\ also \cite{tak:92}, \cite{tak:95}\footnote{The
normalisations in these papers are different. Here we normalise as in
\cite{tak:16}.}) as the gauge transformation matrix $\tilde M$ of the
$L$-matrix.
\comment{note 27 p.3 (4)16/VII/2018}
\begin{equation}
    M_\lambda(v) :=
    \begin{pmatrix}
             - \theta_{00}(\frac{\lambda - v}{2}, \frac{\tau}{2}) &
             - \theta_{00}(\frac{\lambda + v}{2}, \frac{\tau}{2}) \\
    \phantom{-}\theta_{01}(\frac{\lambda - v}{2}, \frac{\tau}{2}) &
    \phantom{-}\theta_{01}(\frac{\lambda + v}{2}, \frac{\tau}{2})
    \end{pmatrix}.
\label{def:gauge-trans}
\end{equation}

\begin{remark}
\label{rem:gauss-decomp-Mlambda}
\comment{note23 p.25 (3)20/X/2017}
 In fact, the lower triangular part of the Gauss decomposition of
 $M_\lambda(v)$ is essentially equal to the transpose of \eqref{tildeM}
 in Baxter's method in \secref{sec:baxter72}. There the conditions
 $p_0=p_1$ and $p_r=p_{r+1}$ require $v$ to be zero, which makes
 $M_\lambda(v)$ degenerate and the Gauss decomposition
 diverges. Nevertheless a part of this Gauss decomposition survives and
 gives the matrix $\tilde M_{i''}$ in \secref{subsec:QR:baxter72},
 \eqref{tildeM}.
\end{remark}

Let us quote several formulae from the previous work \cite{tak:16}. The
action of the elements of the twisted $L$-matrix ((3.2) in
\cite{tak:16}),
\begin{equation}
    L_{\lambda,\lambda'}(u;v) =
    \begin{pmatrix}
    \alpha_{\lambda,\lambda'}(u;v) &  \beta_{\lambda,\lambda'}(u;v) \\ 
    \gamma_{\lambda,\lambda'}(u;v) & \delta_{\lambda,\lambda'}(u;v)
    \end{pmatrix}
    := M_\lambda(v)^{-1} L(u) M_{\lambda'}(v),
\label{def:twisted-L}
\end{equation}
on a vector (a function) in $\Theta^{4l++}_{00}$ (cf.\
\eqref{def:[z;a]}), 
\begin{equation}
    \omega_\lambda(u;v)
    :=
    \left[ z ; \frac{\lambda+u-v}{2}+(-l+1)\eta \right]_{2l},
\label{pseudo-vac}
\end{equation}
is as follows ((3.6) and (3.8) of \cite{tak:16}):
\begin{equation}
\begin{aligned}
    \alpha_{\lambda\pm4l\eta,\lambda}(u;v) \omega_{\pm\lambda}(u;\pm v)
    &= 2 [u + 2l\eta] \omega_{\pm\lambda-2\eta}(u;\pm v),
\\
    \gamma_{\lambda\pm4l\eta,\lambda}(u;v) \omega_{\pm\lambda}(u;\pm v)
    &= 0,
\\
    \delta_{\lambda\pm4l\eta,\lambda}(u;v) \omega_{\pm\lambda}(u;\pm v)
    &=
    2[u-2l\eta]\frac{[\lambda]}{[\lambda\pm4l\eta]}
    \omega_{\pm\lambda+2\eta}(u;\pm v).
\end{aligned}
\label{action-on-vac}
\end{equation}

Note that
\comment{note27 pp.5-6 (1-2)18/VII/2018}
\comment{note27 pp.8-9 (2-3)19/VII/2018}
\begin{gather}
    M_{\lambda+4rl\eta}(v) = M_{\lambda}(v),
\label{Mlambda:period}
\\
    \omega_{\lambda+4rl\eta}(u;v) = \omega_\lambda(u;v),
\label{omega-lambda:period}
\\
    [\lambda+4rl\eta] = [\lambda].
\label{[lambda]:period}
\end{gather}
because $2rl\eta=r'\in\Integer$

Fix parameters $\lambda_0$ and $v$ and define
\begin{equation}
    \tilde M_{i''} := M_{\lambda_0+4i''l\eta}(v).
\label{tildeM:gauge}
\end{equation}
Then, when $i''=j''\pm1\pmod{r}$, the matrix elements $\alpha_{i'',j''}$,
$\gamma_{i'',j''}$, $\delta_{i'',j''}$ of the twisted $L$-matrix in
\eqref{MLM-1} act on $\omega_{\pm\lambda_0\pm4j''l\eta}(u;\pm v)$ as
follows:
\comment{note27 p.4 (1)17/VII/2018}
\begin{equation}
\begin{aligned}
    \alpha_{j''\pm1,j''}(u;v) 
    \omega_{\pm\lambda_0\pm4j''l\eta}(u;\pm v)
    &=
    2 [u + 2l\eta]
    \omega_{\pm\lambda_0\pm4j''l\eta}(u-2\eta;\pm v),
\\
    \gamma_{j''\pm1,j''}(u;v)
    \omega_{\pm\lambda_0\pm4j''l\eta}(u;\pm v)
    &= 0,
\\
    \delta_{j''\pm1,j''}(u;v)
    \omega_{\pm\lambda_0\pm4j''l\eta}(u;\pm v)
    &=
    2[u-2l\eta]
    \frac{[\lambda_0+4j''l\eta]}{[\lambda_0+4(j''\pm1)l\eta]}
    \omega_{\pm\lambda_0\pm4j''l\eta}(u+2\eta;\pm v).
\end{aligned}
\label{twisted-local-L:action}
\end{equation}
(We used $\omega_{\lambda\pm2\eta}(u;v)=\omega_\lambda(u\pm2\eta;v)$.)

Hence, if we define $S(u)=(S^{i''}_{j''}(u))_{i'',j''=1,\dotsc,2r}$ by
\comment{note27 p.3 (4)16/VII/2018}
\begin{equation}
    S^{j''\pm1(\bmod{2r})}_{j''}(u):
    \Comp^r\owns e_k \mapsto
    \tau_{kj''}\, \omega_{\pm\lambda_0\pm4j''l\eta}(u;\pm v)
\label{Si''j'':fab}
\end{equation}
($\{e_k\}_{k=1,\dotsc,2r}$ is the standard basis of $\Comp^{2r}$) and
\begin{equation}
    S^{i''}_{j''}(u) = 0
    \text{ unless } i''=j''\pm1 \pmod{2r},
\label{Sn(u):zero-element:fab}
\end{equation}
(that is to say, $S(u)$ is of the form \eqref{S:matrix:fab}),
we have (cf.\ \eqref{M-1(LS)M})
\begin{equation}
    M^{-1} (L(u)\tensor S(u)) M
    =
    \begin{pmatrix}
     A(u) & * \\ 0 & D(u)
    \end{pmatrix},
\label{M-1(LS)M:fab}
\end{equation}
where $M$ is defined by
\begin{equation}
    M = \sum_{i''=1}^{2r} \tilde M_{i''} \tensor E^{i''}_{i''}
\label{M:fab}
\end{equation}
with $\tilde M_{i''}$ in \eqref{tildeM:gauge}. The matrices $A(u)$
and $D(u)$ in \eqref{M-1(LS)M:fab} are
\comment{note 27 p.4 (1)17/VII/2018}
\begin{equation}
 \begin{aligned}
    A(u)
    &=
    2[u+2l\eta] S(u-2\eta),
\\
    D(u)
    &=
    X^{-1} \bigl(2[u-2l\eta] S(u+2\eta)\bigr) X,
 \end{aligned}
\label{A,D:fab}
\end{equation}
where $X=\diag_{i''=1,\dotsc,r}([\lambda_0+4i''l\eta])$. Therefore, as
in \secref{sec:baxter72}, we obtain the relation \eqref{TQR:temp}, which
reduces to the $TQ$-relation \eqref{TQR} with the same coefficients
$h_\pm(u)$ as before.

\subsection{$Q_L$ and its commutation relation with $Q_R$}
\label{subsec:QL:fabricius}

The construction of the operator $Q_L$ is exactly the same as in
\secref{subsec:QL:baxter72}. Namely, the operator defined by
\eqref{def:QL} satisfies the $QT$-relation \eqref{QLT}.

\comment{note27 p.7 (1)19/VII/2018}

As in \secref{subsec:QLQR=QLQR:baxter72}, we need to show the
commutation relation of $Q_R(u)$ and $Q_L(u')$, \eqref{QLQR=QLQR}. The
strategy of the proof is the same as in
\secref{subsec:QLQR=QLQR:baxter72}: Find a matrix $Y$ satisfying
\eqref{YWY-1} for the matrix $W(j,i|u,v)$ defined by
\eqref{def:W(ji|uv)}. 

The matrix elements $W(j,i|u,v)^{j''',i'''}_{j'',i''}=\langle
S^{j'''}_{j''}(u) e_{j}, S^{i'''}_{i''}(v) e_{i} \rangle$ in the matrix
$W(j,i|u,v)$ are non-zero only when $j'''=j''\pm1\pmod{2r}$ and
$i'''=i''\pm1\pmod{2r}$ and have the following form by the formula
\eqref{<omega,omega>}: 
\comment{note27 pp.10--13, (4)--(7)19/VII/2018}
\begin{equation}
 \begin{split}
    & W(j,i|u,v)^{j''',i'''}_{j'',i''}
\\
    ={}&
    \overline{\tau_{jj''}}\tau_{ii''}C'_{2l}
\\  &{}\times
    \theta^{(2l)}_{00}
    \left( \sqrt{-1}\Im{}+2(i''-j'')l\eta +
           \frac{1}{2}u^{(1)j'''i'''}_{j''i''} +
           \mu^{(1)j'''i'''}_{j''i''}l\eta \right)
\\  &{}\times
    \theta^{(2l)}_{00}
    \left( \Re{}         +2(i''+j'')l\eta +
           \frac{1}{2}u^{(2)j'''i'''}_{j''i''} +
           \mu^{(2)j'''i'''}_{j''i''}l\eta \right).
 \end{split}
\label{W(j,i|u,v):element:fab}
\end{equation}
(See \eqref{def:C',theta2l} for notations $C'_{2l}$ and
$\theta^{(2l)}_{00}$.)  Here, $\Im{}=\Im(\lambda_0-v)$ and
$\Re{}=\Re(\lambda_0-v)$. The parameters $u^{(1)j'''i'''}_{j''i''}$,
$u^{(2)j'''i'''}_{j''i''}$, $\mu^{(1)j'''i'''}_{j''i''}$ and
$\mu^{(2)j'''i'''}_{j''i''}$ are defined in \tableref{tab:u(1),mu(1)}
and \tableref{tab:u(2),mu(2)}.

\begin{table}[h]
\comment{note27 pp.10--13, (4)--(7)19/VII/2018}
\begin{center}
 \begin{tabular}{|l||c|c|}
 \hline
  $u^{(1)j'''i'''}_{j''i''}$, $\mu^{(1)j'''i'''}_{j''i''}$
               & $j'''=j''+1$  & $j'''=j''-1$ \\
 \hline\hline
  $i'''=i''+1$ & $ u'+u$, $ 0$ & $ u'-u$, $ 2$ \\
 \hline
  $i'''=i''-1$ & $-u'+u$, $-2$ & $-u'-u$, $0$ \\
 \hline
 \end{tabular}
\caption{$u^{(1)j'''i'''}_{j''i''}$ and $\mu^{(1)j'''i'''}_{j''i''}$ in
 \eqref{W(j,i|u,v):element:fab}.}
\label{tab:u(1),mu(1)}
\end{center}
\end{table}

\begin{table}[h]
\comment{note27 pp.10--13, (4)--(7)19/VII/2018}
\begin{center}
 \begin{tabular}{|l||c|c|}
 \hline
  $u^{(2)j'''i'''}_{j''i''}$, $\mu^{(2)j'''i'''}_{j''i''}$
               & $j'''=j''+1$  & $j'''=j''-1$ \\
 \hline\hline
  $i'''=i''+1$ & $ u'-u$, $ 2$ & $ u'+u$, $ 0$ \\
 \hline
  $i'''=i''-1$ & $-u'-u$, $ 0$ & $-u'+u$, $-2$ \\
 \hline
 \end{tabular}
\caption{$u^{(2)j'''i'''}_{j''i''}$ and $\mu^{(2)j'''i'''}_{j''i''}$ in
 \eqref{W(j,i|u,v):element:fab}.}
\label{tab:u(2),mu(2)}
\end{center}
\end{table}

To find a diagonal matrix
$Y=(Y^{j'',i''}_{j''',i'''})_{(j'',i''),(j''',i''')\in\{1,\dotsc,r\}^2}$
satisfying 
\begin{equation}
    Y W(j,i|u,u') Y^{-1} = W(j,i|u',u),
\label{YWY-1:fab}
\end{equation}
we have to solve the following four series of equations:
\comment{note27 p.14 (8)19/VII/2018}
\begin{align}
    \frac{y_{i''+1,j''+1}}{y_{i''j''}}
    &=
    \frac{\theta^{(2l)}_{00}(\Re{}+2(i''+j'')l\eta+\frac{u'-u}{2}+2l\eta)}
         {\theta^{(2l)}_{00}(\Re{}+2(i''+j'')l\eta+\frac{u-u'}{2}+2l\eta)}
    =: A_{i''j''},
\label{++}
\\
    \frac{y_{i''+1,j''-1}}{y_{i''j''}}
    &=
    \frac
    {\theta^{(2l)}_{00}(\sqrt{-1}\Im{}+2(i''-j'')l\eta+\frac{u'-u}2+2l\eta)}
    {\theta^{(2l)}_{00}(\sqrt{-1}\Im{}+2(i''-j'')l\eta+\frac{u-u'}2+2l\eta)}
    =: B_{i''j''},
\label{+-}
\\
    \frac{y_{i''-1,j''+1}}{y_{i''j''}}
    &=
    \frac
    {\theta^{(2l)}_{00}(\sqrt{-1}\Im{}+2(i''-j'')l\eta+\frac{-u'+u}2-2l\eta)}
    {\theta^{(2l)}_{00}(\sqrt{-1}\Im{}+2(i''-j'')l\eta+\frac{-u+u'}2-2l\eta)},
\label{-+}
\\
    \frac{y_{i''-1,j''-1}}{y_{i''j''}}
    &=
    \frac{\theta^{(2l)}_{00}(\Re{}+2(i''+j'')l\eta+\frac{-u'+u}2-2l\eta)}
         {\theta^{(2l)}_{00}(\Re{}+2(i''+j'')l\eta+\frac{-u+u'}2-2l\eta)}.
\label{--}
\end{align}
It is easy to show that the equation \eqref{--} is equivalent to
\eqref{++} and the equation \eqref{-+} is equivalent to
\eqref{+-}. Hence what we need to solve is the system \eqref{++} and
\eqref{+-}, namely the system of linear difference equations,
\begin{align}
    y_{i''+1,j''+1} &= A_{i''j''}\, y_{i''j''},
\label{++:linear}
\\
    y_{i''+1,j''-1} &= B_{i''j''}\, y_{i''j''}.
\label{+-:linear}
\end{align}
The compatibility condition of this sytem is
\comment{note27 p.16 (1)20/VII/2018}
\begin{equation}
    B_{i''+1,j''+1}\, A_{i''j''} = A_{i''+1,j''-1}\, B_{i''j''},
\label{++,+-:compatible}
\end{equation}
which is readily checked.

We need solutions of the linear system satisfying the periodic boundary
condition:
\begin{equation}
    y_{i''+2r,j''} = y_{i'',j''+2r} = y_{i''j''}.
\label{y:periodic}
\end{equation}

\begin{lemma}
\label{lem:y:period:diagonal}
 Assume $\lambda_0-v = 2r''l\eta$ ($r''\in \Integer$). Then for any
 $(i'',j'')\in\Integer^2$
\begin{align}
    \prod_{k''=0}^{r-1}A_{i''+k'',j''+k''} &= 1,
\label{prod-A=1}
\\
    \prod_{k''=0}^{r-1}B_{i''+k'',j''-k''} &= 1.
\label{prod-B=1}
\end{align}
\end{lemma}

Hereafter we assume
\begin{equation}
    \lambda_0-v=2r''l\eta \qquad (r''\in\Integer).
\label{lambda0-v=2r''leta}
\end{equation}
Under this assumption \lemref{lem:y:period:diagonal} shows
$y_{i''+r,j''+r}=y_{i''+r,j''-r}=y_{i''j''}$ as a consequence of
\eqref{++} and \eqref{+-}. The periodicity \eqref{y:periodic} directly
follows from this.

Thus the existence of $Y$ satisfying \eqref{YWY-1:fab} has been proved
and therefore the commutation relation \eqref{QLQR=QLQR} has been shown.

\begin{proof}[Proof of \lemref{lem:y:period:diagonal}]
 The proof is similar to the proof of (66) in \cite{fab:07}.
 
 The assumption on $\lambda_0-v$ means
\[
    \Re=\Re(\lambda_0-v)=2r''l\eta\, \qquad
    \Im=\Im(\lambda_0-v) = 0,
\]
 because $\eta\in\Real$.

 By multiplying $A_{i''j''}$'s defined by \eqref{++}, we have
\[
    \prod_{k''=0}^{r-1} A_{i''+k'',j''+k''}
    =
    \prod_{k''=0}^{r-1}
    \frac
    {\theta^{(2l)}_{00}(2(r''+i''+j''+2k''+1)l\eta+\frac{u'-u}{2})}
    {\theta^{(2l)}_{00}(2(r''+i''+j''+2k''+1)l\eta+\frac{u-u'}{2})}.
\]
 Hence the equation \eqref{prod-A=1} holds, if
\begin{multline}
    \{r''+i''+j''+2k''_1+1 \pmod{r} \mid
    k''_1\in\{0,1,\dotsc,r-1\}\}
\\
    =
    \{-(r''+i''+j''+2k''_2+1) \pmod{r} \mid
    k''_2\in\{0,1,\dotsc,r-1\}\},
\label{prod-A=1:condition}
\end{multline}
 because $\theta^{(2l)}_{00}(u)$ is even and has a period
 $2rl\eta\in\Integer$. The condition \eqref{prod-A=1:condition} means
 that for any $k_1''\in\{0,1,\dotsc,r-1\}$ there exists a unique
 $k''_2\in\{0,1,\dotsc,r-1\}$ which satisfies
\[
     r''+i''+j''+2k''_1+1\equiv -r''-i''-j''-2k''_2-1 \pmod{r}.
\]
 It is sufficient to take $k''_2\in\{0,1,\dotsc,r-1\}$ such that
\[
     k''_2 \equiv -(r''+i''+j''+k''_1+1) \pmod{r},
\]
 which exists uniquely. Equation \eqref{prod-A=1} is proved.

 The proof of the equation \eqref{prod-B=1} is similar, because
\[
    \prod_{k''=0}^{r-1} B_{i''+k'',j''-k''}
    =
    \prod_{k''=0}^{r-1}
    \frac
    {\theta^{(2l)}_{00}(2(i''-j''+2k''+1)l\eta+\frac{u'-u}{2})}
    {\theta^{(2l)}_{00}(2(i''-j''+2k''+1)l\eta+\frac{u-u'}{2})},
\]
 thanks to the assumptionn $\Im{}=0$.
\end{proof}

As in \secref{subsec:QLQR=QLQR:baxter72}, if there are non-degenerate
$Q_R(u_0)$ and $Q_L(u_0)$, the $Q$-operator satisfying the $TQ$-relation
\eqref{TQ}, the $QT$-relation \eqref{QT} and commutativity \eqref{QQ=QQ}
can be constructed by \eqref{Q=QRQR-1=QL-1QL}. We do not discuss the
non-degeneracy problem here.

\subsection{Quasi-periodicity of the $Q$-operator}
\label{subsec:quasi-periodicity:fabricius}

The unitary operators $U_a$ ($a=1,3$) act on $\omega_\lambda(u;v)$
defined by \eqref{pseudo-vac} as follows.
\begin{equation}
 \begin{aligned}
\comment{note27 p.22 (1)26/VII/2018}
    U_1 \omega_\lambda(u;v) &= e^{-l \pi i} \omega_\lambda(u+1,v),
\\
\comment{note26 p.4 (2)17/VI/2018}
    U_3 \omega_\lambda(u,v)
    &= e^{l(\tau-1)\pi i +2l(\lambda+u-v+2l\eta)\pi i }
       \omega_\lambda(u+\tau,v).
 \end{aligned}
\label{Ua-on-omega}
\end{equation}
Hence we obtain
\begin{equation}
 \begin{aligned}
    \comment{note27 p.23 (1)27/VII/2018}
    U_1 S^{i''}_{j''}(u)&= e^{-l\pi i} S^{i''}_{j''}(u+1),
\\
    \comment{note27 p.25 (3)27/VII/2018}
    U_3 S^{i''}_{j''}(u)&=
    e^{l(\tau-1)\pi i +2l \pi i u+ 2l\pi i\, d(i'',j'')}
    S^{i''}_{j''}(u+\tau),
 \end{aligned}
\label{UaS(u):fab}
\end{equation}
from \eqref{Si''j'':fab}. Here
\begin{equation}
    d(i'',j'')=
    \begin{cases}
      2(2j''+1)l\eta + (\lambda_0-v) &(i''=j''+1\pmod{2r}), \\
     -2(2i''+1)l\eta - (\lambda_0-v) &(i''=j''-1\pmod{2r}).
     \end{cases}
\label{def:d(i'',j''):fab}
\end{equation}
As in \secref{subsec:quasi-periodicity:baxter72}, the quasi-periodicity
of $Q_R(u)$ with respect to $u\mapsto u+1$ follows directly from
\eqref{UaS(u):fab}:
\begin{equation}
    U_1^{\tensor N} Q_R(u) = e^{-N\pi i l} Q_R(u+1).
\label{U1QR(u)=QR(u+1):fab}
\end{equation}
The quasi-periodicity with respect to $u\mapsto u+\tau$ is derived using
a similarity transformation of $S(u)$ by a diagonal matrix
\comment{note27 p.32 (1)22/X/2018}
\begin{equation}
    S(u) \mapsto A S(u) A^{-1},\quad
    A=\diag_{i''=1,\dotsc,r}(A_{i''}),\
    A_{i''} := \exp(4j''(j''+r'')l^2\eta\pi i),
\label{def:A:fab}
\end{equation}
where $r''$ is the integer in \eqref{lambda0-v=2r''leta}. Note that
$A_{i''+2r}=A_{i''}$ because of the rationality $\eta=\dfrac{r'}{2rl}$.
As a result of \eqref{UaS(u):fab}, we have
\begin{equation}
    U_3^{\tensor N} Q_R(u)
    =
    e^{Nl \pi i (\tau-1)+2 Nl\pi i u} Q_R(u+\tau).
\label{U3QR(u)=QR(u+tau):fab}
\end{equation}
(cf.\ \eqref{U3QR(u)=QR(u+tau):bax}.)

Derivation of the quasi-periodicity of the operator $Q_L$ from that of
$Q_R$ is also the same as in
\secref{subsec:quasi-periodicity:baxter72}. (Note that $Q_R(u+2)=Q_R(u)$
follows from $\omega_\lambda(u+2,v)=\omega_\lambda(u,v)$.)

Thus we obtain the quasi-periodicity of the $Q$-operator
\begin{equation}
 \begin{aligned}
    U_1^{\tensor N} Q(u) &= 
    Q(u) U_1^{\tensor N}
    = e^{-Nl\pi i} Q(u+1),
\\
    U_3^{\tensor N} Q(u) &= 
    Q(u) U_3^{\tensor N}
    = e^{Nl\pi i(\tau-1)+2Nl\pi iu} Q(u+\tau).
 \end{aligned}
\label{Q:quasi-periodicity:fab}
\end{equation}

\section{Concluding remarks}
\label{sec:conclusion}

We have constructed the $Q$-operator of the higher spin eight vertex
model satisfying
\begin{itemize}
 \item $TQ$- and $QT$-relations, \eqref{TQ}, \eqref{QT},
 \item commutativity \eqref{QQ=QQ},
 \item quasi-periodicity, \eqref{Q:quasi-periodicity:bax} and
       \eqref{Q:quasi-periodicity:fab},
\end{itemize}
in two ways. The computations are far more complicated, but,
surprisingly the strategies by Baxter and Fabricius for the spin $1/2$
case (the eight vertex model) work also in higher spin cases almost as
they are. This fact strongly suggests that a representation theoretical
structure is hidden behind those technical complicated construction.

One clue could be the connection of Fabricius's construction and the
generalised algebraic Bethe Ansatz mentioned in
\remref{rem:gauss-decomp-Mlambda}. 

Let us make several comments related to other works:
\begin{itemize}
 \item As in our previous work \cite{tak:16} the Bethe Ansatz equation
       for the eigenvalues of the transfer matrix and the sum rule of
       the Bethe roots (integrality of the sum of the Bethe roots)
       follow from the above mentioned properties of the $Q$-operator.
       We do not repeat the same derivation here.

 \item The $Q$-operators for the XXZ spin chain of higher spin were
       constructed by Roan in \cite{roa:07-1}, following both ways of
       Baxter, \cite{bax:72} and \cite{bax:73}. The first method is
       similar to the method in \secref{sec:fabricius}. In the same
       paper other functional relation for the $Q$-operator and fused
       transfer matrices are derived, which should be generalised to our
       elliptic model. (See also \cite{mot:13} for a related work.)

 \item Roan's construction in \cite{roa:07-2} of the $Q$-operator for
       the eight vertex model is similar to that by Fabricius
       \cite{fab:07}. The main difference is the valudes of
       parameters. In this context our construction in this paper adds
       another parameter $\lambda_0$ by the use of the matrix
       $M_\lambda(v)$. 
      
 \item The $Q$-operators for the elliptic models with
       infinite-dimensional state spaces were constructed by Zabrodin in
       \cite{zab:00} and by Chicherin, Derkachov, Karakhanyan and
       Kirschner in \cite{c-d-k-k:13}. The latter construction
       (especially in \S2 of that paper) seems to have something in
       common with ours. The relation of our construction and reduction
       of those operators in \cite{zab:00}, \cite{c-d-k-k:13},
       \cite{c-d-s:14} would be important.
       
 \item The non-locality problem mentioned in the last section of
       \cite{tak:16} is also present in the construction of this paper,
       which is inevitable, if one takes Baxter's strategy. It is an
       interesting question whether completely
       different approach as in \cite{b-l-z:97}, \cite{b-l-z:99},
       \cite{man:14-1}, \cite{man:14-2}, \cite{bazh-str:90} would be
       applicable. 

\end{itemize}

\subsection*{Acknowledgements}

The author expresses his gratitude to Hitoshi Konno, Kohei Motegi for
discussions and encouragement and to Klaus Fabricius for informing
references.

The author is grateful to Rikkyo University and Tokyo Unviersity
of Marine Science and Technology for there hospitality, where parts of
this work were done.

This work has been funded by the Russian Academic Excellence Project
`5--100'.

\begin{appendices}
\section{Sklyanin algebra}
\label{app:sklyanin}

In this appendix we recall several facts on the Sklyanin algebra and
its representations from \cite{skl:82} and \cite{skl:83}.

The {\em Sklyanin algebra} is an associative algebra generated by four
generators $S^a$ ($a=0,\dotsc,3$) subject to the following relations:
\begin{equation}
    L_{12}(v) L_{13}(u) R_{23}(u-v) = R_{23}(u-v) L_{13}(u) L_{12}(v).
\label{RLL}
\end{equation}
Here the symbols are defined as follows:
\begin{itemize}
 \item the {\em $L$-operator} $L(u)$ with a complex parameter $u$ is
defined by
\begin{equation}
    L(u) = \sum_{a=0}^3 W_a^L(u) S^a \tensor \sigma^a,
\label{def:L:alg}
\end{equation}
where
\begin{equation}
\begin{aligned}
    W_0^L(u) 
    &= \frac{\theta_{11}(u,\tau)}{\theta_{11}(\eta,\tau)},&
    W_1^L(u)
    &= \frac{\theta_{10}(u,\tau)}{\theta_{10}(\eta,\tau)},
\\
    W_2^L(u)
    &= \frac{\theta_{00}(u,\tau)}{\theta_{00}(\eta,\tau)},&
    W_3^L(u)
    &= \frac{\theta_{01}(u,\tau)}{\theta_{01}(\eta,\tau)}.
\end{aligned}
\label{def:Wa}
\end{equation}
 \item The matrix $R(u)$ is {\em Baxter's $R$-matrix} defined by
\begin{equation}
    R(u) = \sum_{a=0}^3 W_a^R(u) \sigma^a \tensor \sigma^a,\qquad
    W_a^R(u) := W_a^L(u + \eta).
\label{def:R}
\end{equation}
       Explicitly, it has the form (cf.\ \cite{tak:95} Appendix A)
\begin{equation}
    R(u) = \begin{pmatrix}
     a(u) &    0 &    0 & d(u)\\
        0 & b(u) & c(u) & 0   \\
        0 & c(u) & b(u) & 0   \\
     d(u) &    0 &    0 & a(u)
     \end{pmatrix},
\label{R}
\end{equation}
       where 
\begin{align*}
    a(u) &=
    C \theta_{01}(2it\eta,2it)\,\theta_{01}(itu,2it)
    \,\theta_{11}(it(u+2\eta),2it) \\
    b(u) &=
    C \theta_{01}(2it\eta,2it)\,\theta_{11}(itu,2it)
    \,\theta_{01}(it(u+2\eta),2it) \\
    c(u) &=
    C \theta_{11}(2it\eta,2it)\,\theta_{01}(itu,2it)
    \,\theta_{01}(it(u+2\eta),2it) \\
    d(u) &=
    C \theta_{11}(2it\eta,2it)\,\theta_{11}(itu,2it)
    \,\theta_{11}(it(u+2\eta),2it) \\
    C &= \frac
    {-2 e^{-\pi t u(u+2\eta)}}
    {\theta_{01}(0,2it)\,
     \theta_{01}(2it\eta,2it)\,\theta_{11}(2it\eta,2it)},\quad
    t=\frac{i}{\tau}.
\end{align*}
 \item The indices designate the spaces on which operators act
       non-trivially. For example,
\begin{equation*}
    L_{12}(u) =
    \sum_{a=0}^3
    W_a^L(u) S^a \tensor \sigma^a \tensor 1,\qquad
    R_{23}(u) = 
    \sum_{a=0}^3
    W_a^R(u) 1 \tensor \sigma^a \tensor \sigma^a.
\end{equation*}
\end{itemize} 

Although the relation \eqref{RLL} contains parameters $u$ and $v$, the
commutation relations among $S^a$ ($a=0, \dots, 3$) do not depend on
them:
\begin{equation}
    [S^\alpha, S^0    ]_- =
    -i J_{\alpha,\beta} [S^\beta,S^\gamma]_+, \qquad
    [S^\alpha, S^\beta]_- =
                      i [S^0,    S^\gamma]_+,
\label{comm_rel}
\end{equation}
where $(\alpha, \beta, \gamma)$ stands for an arbitrary cyclic
permutation of $(1,2,3)$ and $[A,B]_\pm$ are the (anti-)commutator
$AB\pm BA$. The structure constants
$
    J_{\alpha,\beta}
    =
    ((W^L_\alpha)^2-(W^L_\beta)^2)/((W^L_\gamma)^2-(W^L_0)^2)
$
depend on $\tau$ and $\eta$ but not on $u$.

Let $l$ be a positive half integer. The {\em spin $l$ representation}
$\rho^{l}$ of the Sklyanin algebra is defined as follows: The
representation space is a space of entire functions,
\begin{multline}
    \Theta^{4l+}_{00} :=
    \{f(z) \, |\, \\
     f(z+1) = f(-z) = f(z), f(z+\tau)=\exp^{-4l\pi i(2z+\tau)}f(z) \},
\label{def:theta-space}
\end{multline}
which is of dimension $2l+1$. The generator $S^a$ of the Sklyanin
algebra acts as a difference operator on this space:
\begin{equation}
    (\rho^l(S^a) f)(z) =
    \frac{s_a(z-l\eta)f(z+\eta)-s_a(-z-l\eta)f(z-\eta)}
         {\theta_{11}(2z,\tau)},
\end{equation}
where
\begin{alignat*}{2}
    s_0(z) &=  \theta_{11}(\eta,\tau) \theta_{11}(2z,\tau),\qquad&
    s_1(z) &=  \theta_{10}(\eta,\tau) \theta_{10}(2z,\tau),\\
    s_2(z) &= i\theta_{00}(\eta,\tau) \theta_{00}(2z,\tau),\qquad&
    s_3(z) &=  \theta_{01}(\eta,\tau) \theta_{01}(2z,\tau).
\end{alignat*}
In the simplest case $l= 1/2$, $\rho^{1/2}(S^a)$ are expressed by the
Pauli matrices $\sigma^a$. We can identify $\Theta^{2+}_{00}$ and
$\Comp^2$ by
\begin{equation}
\begin{aligned}
    \theta_{00}(2z,2\tau)-\theta_{10}(2z,2\tau) 
    &\longleftrightarrow
    \begin{pmatrix} 1 \\ 0 \end{pmatrix},
\\
    \theta_{00}(2z,2\tau)+\theta_{10}(2z,2\tau) 
    &\longleftrightarrow
    \begin{pmatrix} 0 \\ 1 \end{pmatrix}.
\end{aligned}
\label{rep:identify}
\end{equation}
Under this identification $S^a$ have matrix forms
\begin{equation}
    \rho^{1/2}(S^a) = \theta_{11}(2\eta,\tau) \sigma^a.
\label{rep:pauli}
\end{equation}

The representation space $\Theta^{4l+}_{00}$ has a natural Hermitian
structure defined by the following {\em Sklyanin form}:
\begin{equation}
    \left< f(z), g(z)\right>
    :=
    \int_0^1 dx \int_0^{\tau/i} dy\,
    \overline{f(z)} g(z) \mu(z,\bar z),
\label{skl-form}
\end{equation}
where $z=x+iy$ and the kernel function $\mu(z,w)$ is defined by
\begin{equation}
    \mu(z,w)
    :=
    \frac
    {\theta_{11}(2z,\tau) \theta_{11}(2w,\tau)}
    {\displaystyle
     \prod_{j=0}^{2l+1}
     \theta_{00}(z+w+(2j-2l-1)\eta,\tau)\theta_{00}(z-w+(2j-2l-1)\eta,\tau)
    }.
\label{def:mu}
\end{equation}
The most important property of this sesquilinear positive definite 
scalar product is that the generators $S^a$ of the Sklyanin algebra
become self-adjoint:
\begin{equation}
    (S^a)^* = S^a,\text{ namely, }
    \langle f(z), S^a g(z) \rangle = \langle S^a f(z), g(z) \rangle.
\label{Sa:self-adj}
\end{equation}

In \cite{skl:83} Sklyanin also defined involutive automorphisms:
\begin{equation}
    X_a: (S^0, S^a,  S^b,  S^c) \mapsto
         (S^0, S^a, -S^b, -S^c),
\label{def:X-a}
\end{equation}
for $a=1,2,3$, where $(a,b,c)$ is a cyclic permutation of $(1,2,3)$. The
unitary operators $U_a$ defined by
\begin{equation}
\begin{aligned}
    U_1: &\Theta^{4\ell+}_{00} \owns f(z) \mapsto &
         &(U_1 f)(z) = e^{\pi i \ell} f\left(z + \frac{1}{2} \right),
\\
    U_3: &\Theta^{4\ell+}_{00} \owns f(z) \mapsto &
         &(U_3 f)(z) = e^{\pi i \ell} e^{\pi i \ell (4z+\tau)}
                       f\left(z + \frac{\tau}{2} \right),
\end{aligned}
\label{def:Ua}
\end{equation}
and $U_2= U_3 U_1$, intertwine representations $\rho^\ell \circ X_a$ and
$\rho^{\ell}$: $\rho^\ell(X_a(S^b)) = U_a^{-1} \rho^\ell(S^b)
U_a$. Operators $U_a$ satisfy the relations: $U_a^2 = (-1)^{2\ell}$,
$U_a U_b = (-1)^{2\ell} U_b U_a = U_c$.

\section{Values of Sklyanin forms of elements of $\Theta^{4l++}_{00}$}
\label{app:sklyanin-form}

In Appendix B of \cite{tak:16} (equation (B.14)), we have computed the
Sklyanin form of two shifted products of theta functions, using the
results by Rosengren, \cite{ros:04} and \cite{ros:07} (see also Konno's
work \cite{kon:05}): 
\begin{equation*}
 \begin{split}
    &\langle
     [z;\alpha]_N, [z;\gamma]_N
    \rangle
\\
    ={}&
    C_{N} e^{\pi i N \tau/2}
    \prod_{j=0}^{N-1}
    \theta_{00}(\gamma-\bar\alpha + (2j-N+1)\eta,\tau)
    \theta_{00}(\gamma+\bar\alpha + (2j+N-1)\eta,\tau),
 \end{split}
\end{equation*}
where
\begin{equation*}
    C_N
    =
    \frac{-2\eta e^{3\pi i \tau/4}}
         {[2(N+1)\eta] \prod_{j=1}^\infty (1-e^{2j \pi i \tau})^3}.
\end{equation*}
Hence, the Sklyanin form of $f_\epsilon(\lambda,-\bar u,z)$ and
$f_{\epsilon'}(\lambda',v,z)$ ($\epsilon,\epsilon'=\pm$,
$\lambda,\lambda'\in\Real$, cf.\ \eqref{null-gamma(lambda)}) has the
following form:
\begin{equation}
 \begin{split}
    &\langle
     f_\epsilon(\lambda,-\bar u,z), f_{\epsilon'}(\lambda',v,z)
    \rangle
\\
    ={}&
    F\left(
      \frac{\epsilon'\lambda-\epsilon\lambda}{2}+\frac{v+u}{2}
     \right)
    G\left(
      \frac{\epsilon'\lambda+\epsilon\lambda}{2}+\frac{v-u}{2}
      + 2(2l-1)\eta
     \right),
 \end{split}
\label{<f,f>}
\end{equation}
where $F$ and $G$ are defined by
\begin{align}
    F(z) &:=
    C_{2l} e^{\pi i l \tau}
    \prod_{j=0}^{2l-1} \theta_{00}(z + (2j-2l+1)\eta,\tau),
\label{def:F}
\\
    G(z) &:=
    \prod_{j=0}^{2l-1}
    \theta_{00}(z + (2j+2l-1)\eta - 2(2l-1)\eta,\tau).
\label{def:G}
\end{align}
We shifted the argument in $G(z)$ so that $G(z)$ becomes an even
function:
\begin{equation}
    G(-z) = G(z).
\label{G:even}
\end{equation}
It has also the periodicity:
\begin{equation}
    G(z+1)=G(z),
\label{G:periodic}
\end{equation}
because of the periodicity of $\theta_{00}$.

\medskip
In \secref{subsec:QL:fabricius} we need the Sklyanin form among
$\omega_\lambda(u,v)$'s defined by \eqref{pseudo-vac}. The following
formula is useful.
\comment{note27 p.9 (3) 19/VII/2018}
\begin{equation}
 \begin{split}
    &\langle
     \omega_{\sigma \lambda }(-\bar u,\sigma  v ),
     \omega_{\sigma'\lambda'}(     u',\sigma' v')
    \rangle
\\
    ={}&
    C'_{2l}\,
    \theta^{(2l)}_{00}
    \left(\frac{(\lambda'-v')-\overline{(\lambda-v)}}{2}+
          \frac{\sigma'u'+\sigma u}{2}+(\sigma'-\sigma)l\eta \right)
\\
    &\times
    \theta^{(2l)}_{00}
    \left(\frac{(\lambda'-v')+\overline{(\lambda-v)}}{2}+
          \frac{\sigma'u'-\sigma u}{2}+(\sigma'+\sigma)l\eta \right),
 \end{split}
\label{<omega,omega>}
\end{equation}
where
\begin{equation}
 \begin{aligned}
    C'_{2l}&:=C_{2l} e^{\pi i l \tau},
\\
    \theta^{(2l)}_{00}(u)&:=
    \prod_{j=0}^{2l-1}
    \theta_{00}( u + (2j-2l+1)\eta, \tau ).
 \end{aligned}
\label{def:C',theta2l}
\end{equation}

\end{appendices}



\begin{thebibliography}{Bax2}
\bibitem[B1]{bax:72}
Baxter, R. J.:
Partition Function of the Eight-Vertex Lattice Model,
{\em Ann. Phys.} {\bf 70}
(1972),
193--228.

\bibitem[B2]{bax:73} 
Baxter, R. J.:
Eight-Vertex Model in Lattice Statistics and One-Dimensional
Anisotropic Heisenberg Chain I, II, III, 
{\em Ann. Phys.} {\bf 76}
(1973),
1--24, 25--47, 48--71.

\bibitem[B3]{bax:82}
Baxter, R. J.:
Exactly solved models in statistical mechanics.
{\em Academic Press, Inc., London}, 
(1982),
xii+486 pp.

\bibitem[BLZ1]{b-l-z:97}
Bazhanov, V. V., Lukyanov, S. L. and Zamolodchikov, A. B.:
Integrable structure of conformal field theory II, 
$Q$-operator and DDV equation,
{\em Comm. Math. Phys.} {\bf 190}
247--278
(1997)

\bibitem[BLZ2]{b-l-z:99}
Bazhanov, V. V., Lukyanov, S. L. and Zamolodchikov, A. B.:
Integrable structure of conformal field theory III,
The Yang-Baxter relation,
{\em Comm. Math. Phys.} {\bf 200}
297--324
(1999)

\bibitem[BS]{bazh-str:90}
Bazhanov, V. V. and Stroganov, Yu. G.:
Chiral Potts model as a descendant of the six-vertex model,
{\em J. Statist. Phys.} {\bf 59}
799--817
(1990)

\bibitem[CDKK]{c-d-k-k:13}
Chicherin, D., Derkachov, S., Karakhanyan, D. and Kirschner, R.:
Baxter operators with deformed symmetry,
{\em Nuclear Phys.} {\bf B 868} 
(2013),
652–683. 

\bibitem[CDS]{c-d-s:14}
Chicherin, D., Derkachov, S. E. and Spiridonov, V. P.:
New elliptic solutions of the Yang-Baxter equation,
{\em Comm. Math. Phys.} {\bf 345}
(2016), 
507–543.
	     
\bibitem[F]{fab:07}
Fabricius, K.:
A new Q-matrix in the eight-vertex model,
{\em J. Phys.} {\bf A 40}
4075--4086,
(2007).

\bibitem[FH]{fre-her:15}
Frenkel, E., Hernandez, D.:
Baxter's relations and spectra of quantum integrable models,
{\em Duke Math. J.} {\bf 164},
2407--2460,
(2015).

\bibitem[FM1]{fab-mcc:03-05}
Fabricius, K. and McCoy, B. M.:
New developments in the eight vertex model,
{\em J. Statist. Phys.} {\bf 111},
323--337
(2003);
ditto II, Chains of odd length,
{\em J. Stat. Phys.} {\bf 120}
37--70,
(2005).

\bibitem[FM2]{fab-mcc:07}
Fabricius, K. and McCoy, B. M.:
The TQ equation of the eight-vertex model for complex elliptic roots of
unity,
{\em J. Phys.} {\bf A 40},
14893--14926,
(2007).

\bibitem[FM3]{fab-mcc:09}
Fabricius, K. and McCoy, B. M.:
New Q matrices and their functional equations for the eight
vertex model at elliptic roots of unity. 
{\em J. Stat. Phys.} {\bf 134},
643–668,
(2009).


\bibitem[K]{kon:05}
Konno, H.:
The vertex-face correspondence and the elliptic $6j$-symbols,
{\em Lett. Math. Phys.} {\bf 72}
(2005),
243--258. 

\bibitem[Ma1]{man:14-1}
Mangazeev, V. V.:
On the Yang-Baxter equation for the six-vertex model,
{\em Nuclear Phys.} {\bf B 882}
70--96
(2014)

\bibitem[Ma2]{man:14-2}
Mangazeev, V. V.:
$Q$-operators in the six-vertex model,
{\em Nuclear Phys.} {\bf B 886}
166--184
(2014)

\bibitem[Mo]{mot:13}
Motegi, K.:
On Baxter's Q operator of the higher spin XXZ chain at the
Razumov-Stroganov point,
{\em J. Math. Phys.} {\bf 54}
063510, 13 pp. 
(2013)

\bibitem[Mu]{mum:83}
Mumford, D.:
Tata Lectures on Theta I,
{\em Progress in Mathematics} {\bf 28},
Birkhäuser Boston, Inc.,
(1983),
xiii+235 pp. 

\bibitem[Roa1]{roa:07-1}
Roan, S.-S.:
On $Q$-operators of XXZ Spin Chain of Higher Spin,
{\tt arXiv:cond-mat/0702271}.

\bibitem[Roa2]{roa:07-2}
Roan, S.-S.:
The $Q$-operator and functional relations of the eight-vertex model at
root-of-unity $\eta=\frac{2mK}{N}$ for odd $N$,
{\em J. Phys.} {\bf A 40},
11019--11044,
(2007).

\bibitem[Ros1]{ros:04}
Rosengren, H.:
Sklyanin invariant integration,
{\em Int. Math. Res. Not.} {\bf  60} 
(2004), 
3207--3232.

\bibitem[Ros2]{ros:07}
Rosengren, H.:
An elementary approach to 6j-symbols (classical, quantum, rational,
trigonometric, and elliptic),
{\em Ramanujan J.} {\bf 13}
(2007),
131--166. 

\bibitem[S1]{skl:82}
Sklyanin, E. K.:
Some Algebraic Structures Connected with the Yang-Baxter Equation,
{\em Funkts. analiz i ego Prilozh.} {\bf 16-4}
(1982),
27--34,
(in Russian);
{\em Funct. Anal. Appl.} {\bf 16}
(1983),
263--270
(English transl.).

\bibitem[S2]{skl:83}
Sklyanin, E. K.:
Some Algebraic Structures Connected with the Yang-Baxter Equation.
Representations of Quantum Algebras,
{\em Funkts. analiz i ego Prilozh.} {\bf 17-4}
(1983),
34--48,
(in Russian);
{\em Funct. Anal. Appl.} {\bf 17}
(1984),
273--284,
(English transl.).

\bibitem[T1]{tak:92}
Takebe, T.:
Generalized Bethe Ansatz with the general spin representations of the
Sklyanin algebra,
{\em J. Phys.} {\bf A 25}
(1992), 
1071--1083. 

\bibitem[T2]{tak:95}
Takebe, T.:
Bethe ansatz for higher spin eight-vertex models,
{\em J. Phys.} {\bf A 28}
(1995),
6675--6706;
Corrigendum
{\em J. Phys.} {\bf A 29}
(1996),
1563--1566. 

\bibitem[T3]{tak:96}
Takebe, T.:
Bethe ansatz for higher-spin XYZ models --- low-lying excitations, 
{\em J. Phys.} {\bf A 29}
(1996),
6961--6966. 

\bibitem[T4]{tak:16}
Takebe, T.:
$Q$-operators for higher spin eight vertex models with an even number of
sites,
{\em Lett. Math. Phys.} {\bf 106}
(2016),
319–340. 

\bibitem[TF]{takh-fad:79}
Takhtajan, L. A. and Faddeev, L. D.:
The quantum method of the inverse problem and the Heisenberg XYZ model,
{\em Uspekhi Mat. Nauk} {\bf 34:5}
(1979),
13--63
(in Russian);
{\em Russian Math. Surveys} {\bf 34:5}
(1979),
11--68,
(English translation).

\bibitem[WW]{whi-wat}
Whittaker, E. T. and Watson, G. N.:
A course of modern analysis, 
An introduction to the general theory of infinite processes and 
of analytic functions; with an account of the principal transcendental
functions,
the fourth edition,
{\em Cambridge University Press},
(1927),
vi+608pp.

\bibitem[Z]{zab:00}
Zabrodin, A.:
Commuting difference operators with elliptic coefficients from Baxter's
vacuum vectors,
{\em  J. Phys.} {\bf A 33}
(2000),
3825--3850. 

\end{thebibliography}
\end{document}